\documentclass[12pt,reqno]{amsart}
\usepackage{amsaddr}
\usepackage{amssymb,epsf}
\usepackage{amstext,amsmath,amsfonts,amscd,amsthm,graphicx}
\usepackage{upref}
\usepackage{epsfig}
\usepackage[utf8]{inputenc}
\usepackage{latexsym}
\usepackage{eucal}
\usepackage{wrapfig}
\usepackage{color}
\usepackage{multirow}

\newfont{\fnt}{cmsy10}
\newfont{\sss}{cmti10}

\theoremstyle{definition}
\newtheorem{rmk}{Remark}[section]

\title[Phase-Plane vs. Trend-Based Qualitative Analysis in IS--LM]{Comparison of Phase-Plane and Trend-Based Qualitative Analysis of a Dynamic Two-Dimensional IS--LM Model}

\author{Barbora~Voln\'{a}$^1$*, Mirko~Dohnal$^2$}

\address{$^1$Silesian University in Opava, Mathematical Institute in Opava, \\ Na Rybn\'{i}\v{c}ku 1, 746 01 Opava, Czech Republic.
\\ $^2$Brno University of Technology, Faculty of Business and Management, \\ Kolejn\'{i} 2906/4, 61200 Brno, Czech Republic.
\\ *Corresponding author}
\email{Barbora.Volna@math.slu.cz}
\email{dohnal@vutbr.cz}

\keywords{phase-plane analysis, trend-based analysis, qualitative dynamics, IS--LM model}
\subjclass[2020]{34C60, 91B64. \\ \indent {\em JEL Classification.} C02, C63, E12}

\begin{document}

\begin{abstract}
We compare the classical phase-plane analysis of a dynamic two-dimensional \mbox{IS--LM} model with a sign-based trend approach. Applying both methods to the same specification, we examine which qualitative properties of the system can be established solely from sign relations and how these methods can complement each other. In summary, these two approaches capture complementary perspectives on the qualitative behaviour of the considered model, combining the detailed local classification of equilibrium dynamics provided by phase-plane analysis with the broader admissible qualitative structures captured by trend-based analysis.
\end{abstract}

\maketitle

\section{Introduction}

In this paper, we examine and compare two approaches to analysing the dynamical behaviour of a two-dimensional \mbox{IS--LM} model. The first approach is grounded in phase-plane analysis, a well-established framework for studying qualitative properties of planar dynamical systems through linearization and spectral analysis of the Jacobian matrix. The second approach relies solely on sign-based trend analysis, which does not employ eigenvalue classification and operates under substantially weaker informational assumptions. We apply both approaches to a dynamic two-dimensional \mbox{IS--LM} model, a well-known and extensively studied Keynesian macroeconomic model, thereby illustrating their use and allowing for a transparent comparison.

The original \mbox{IS--LM} model, introduced by \cite{hicks}, describes aggregate macroeconomic equilibrium (or disequilibrium) as the simultaneous equilibrium of the goods and money markets. In its dynamic formulation, see, e.g., \cite{gandolfo} and \cite{torre}, the model is typically analysed with respect to system stability. Although the original structure is two-dimensional, numerous extensions to three- or higher-dimensional settings have been developed, incorporating, for instance, government budget constraints, capital accumulation dynamics, or inflation; see, e.g., \cite{barakova}, \cite{bella_mattana_venturi}, \cite{cai}, \cite{decesare_sportelli}, \cite{fanti_manfredi}, \cite{guirao_garcia-rubio_vera}, \cite{neri_venturi}, \cite{rajpal_bhatia_kumar}, \cite{riad_hattaf_yousfi}, \cite{zhou_li_2008}, and \cite{zhou_li_2009}. Here we focus on the two-dimensional dynamic \mbox{IS--LM} framework based on the original model. Higher-dimensional extensions are therefore not considered. Furthermore, we distinguish between linear, general nonlinear, and Kaldor-type formulations. The model is introduced and represented in both its analytical and trend-based forms, i.e., via a system of differential equations and a sign-based representation, see Section~\ref{the_model}.

In this paper, the model is analysed qualitatively. To minimise possible misunderstandings, the term qualitative analysis therefore refers either to phase-plane analysis based on linearization and eigenvalue classification or to the trend-based approach using trend quantifiers (increasing, constant, decreasing).

The two-dimensional dynamic \mbox{IS--LM} system has been extensively analysed using phase-plane and related techniques. Linearization around stationary points, analysis of the Jacobian matrix, and classification of singularities constitute the standard analytical toolkit for investigating local stability and qualitative behaviour. These methods provide detailed information about the system dynamics and have led to a well-established understanding of its equilibrium properties. A concise overview of the phase-plane results relevant to the model is presented in Section~\ref{phase_plane_analysis}.

In contrast, considerably less attention has been devoted to analytical approaches that rely solely on trend information, such as the monotonicity and curvature of the underlying functions. Such reasoning deliberately refrains from spectral analysis of the Jacobian matrix and does not employ eigenvalue-based classification. Instead, it operates under substantially weaker informational assumptions. The roots of a modern qualitative analysis of the trends can already be traced back to \cite{lopmann_bogen}, and later to \cite{hayes} and \cite{dekleer_brown}. Motivated by these ideas, this type of qualitative analysis, understood as an artificial intelligence method, has been developed over the past decades as a set of algorithms relying on sign information that encodes the monotonicity and curvature of the underlying functions; see, e.g., \cite{dohnal_88}, \cite{dohnal_91}, \cite{dohnal_92}, \cite{doubravsky_dohnal}, \cite{hurme_dohnal}, \cite{kesh_raja}, \cite{konecny_vicha_dohnal}, \cite{parsons}, \cite{parsons_dohnal}, \cite{trave-massuyes_ironi_dague}, or \cite{vicha_dohnal}. We apply this perspective to the dynamic two-dimensional \mbox{IS--LM} model in its linear, general nonlinear, and Kaldor-type formulations and derive new results, including scenario sets and transition graphs, see Section~\ref{trend_based_qualitative_analysis}.

Finally, we compare these two qualitative analyses of the dynamic two-dimensional \mbox{IS--LM} model and discuss their implications in Section~\ref{comparative_analysis}. The aim is to examine which qualitative properties can be established solely from sign information and how phase-plane and trend-based analyses provide complementary perspectives on the dynamics of the considered model.

\section{The Model}
\label{the_model}
In this section, the two-dimensional dynamic \mbox{IS--LM} model is presented in its classical analytical form, i.e., as a system of differential equations, and in its sign-based representation.

\subsection{The Model via System of Differential Equations}
\label{system_of_differential_equation}

A general dynamic two-dimensional \mbox{IS--LM} model, see, e.g., \cite{gandolfo}, is given by the system of differential equations~\eqref{IS-LM_model}. The model describes the simultaneous equilibrium in the goods market (IS~side) and the money market (LM~side).
\begin{equation}
\label{IS-LM_model}
\begin{array}{llcl}
\text{IS:} & \frac{dy}{dt} &=& \alpha \bigl[i(y,r) - s(y,r)\bigr] \\
\text{LM:} & \frac{dr}{dt} &=& \beta \bigl[l(y,r)-m \bigr]
\end{array}
\end{equation}
where $t$ denotes time, $\alpha,\beta>0$, and $y>0$ and $r \in \mathbb{R}$ are the model variables representing aggregate income (GDP, GNP) and the interest rate, respectively. The assumption of positive $y$ follows from the definition of aggregate income. The functions $i(y,r)$ and $s(y,r)$ denote the investment and saving functions associated with the goods market, while $l(y,r)$ denotes the money demand function and $m>0$ the money stock associated with the money market.

The usual economic properties of the functions $i(y,r)$, $s(y,r)$, and $l(y,r)$, see, e.g., \cite{gandolfo}, are expressed by the following inequalities:
\begin{equation}
\label{economic_is_y} 
  \frac{\partial i}{\partial y}>0, \frac{\partial s}{\partial y}>0,
\end{equation}
\begin{equation}
\label{economic_is_r} 
\frac{\partial i}{\partial r}<0, \frac{\partial s}{\partial r}>0,
\end{equation}
\begin{equation}
\label{economic_l_yr} 
\frac{\partial l}{\partial y}>0, \frac{\partial l}{\partial r}<0.
\end{equation}

One of the simplest ways to incorporate the economic conditions \eqref{economic_is_y}, \eqref{economic_is_r} and \eqref{economic_l_yr} is to assume linear relations of the following form:
\begin{equation}
\label{linear_relations}
\begin{aligned}
i(y,r) &= a\,y - b\,r + i_a, \\
s(y,r) &= \sigma\,y + d\,r + s_a, \\
l(y,r) &= k\,y - h\,r,
\end{aligned}
\end{equation}
where $a,b>0$ are the sensitivity coefficients of investment, $i_a>0$ denotes autonomous investment, $\sigma,d>0$ are the sensitivity coefficients of saving, $-s_a>0$ (i.e., $s_a<0$) represents autonomous saving, and $k,h>0$ are the sensitivity coefficients of money demand. Substituting the linear relations \eqref{linear_relations} into \eqref{IS-LM_model} yields the linear IS--LM model \eqref{linear_IS-LM_model}.
\begin{equation}
\label{linear_IS-LM_model}
\begin{array}{llcl}
\text{IS:} & \frac{d y}{d t} &=& \alpha \bigl[ (a-\sigma)y -(b+d)r +(i_a-s_a) \bigr] \\
\text{LM:} & \frac{d r}{d t} &=& \beta \bigl[ k y-h r-m \bigr] 
\end{array}
\end{equation} 

\begin{figure}[ht]
  \centering
  \includegraphics[height=6cm]{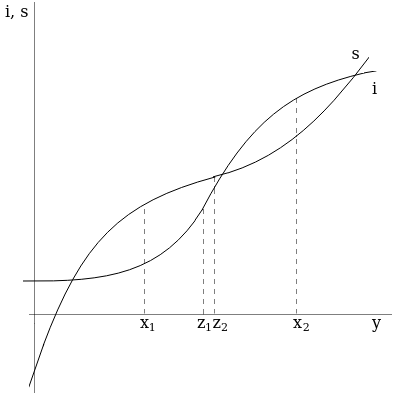}
  \caption{Illustration of the Kaldor conditions}
  \label{fig:Kaldor_is}
\end{figure}

The nonlinear relationships of the model functions can be described using only the economic properties \eqref{economic_is_y}, \eqref{economic_is_r}, and \eqref{economic_l_yr}. 
In addition, we consider the Kaldor-type conditions \cite{kaldor}, defined by the inequalities \eqref{Kaldor_is}, in conjunction with \eqref{economic_is_y}, for the investment and saving functions with respect to aggregate income. These conditions are illustrated in Figure~\ref{fig:Kaldor_is}.
\begin{equation}
\label{Kaldor_is}
  \begin{array}{lll}
  \frac{\partial i}{\partial y} < \frac{\partial s}{\partial y} & \textrm{for} & y \in (-\infty,x_1), \\
  \frac{\partial i}{\partial y} > \frac{\partial s}{\partial y} & \textrm{for} & y \in (x_1,x_2), \\
  \frac{\partial i}{\partial y} < \frac{\partial s}{\partial y} & \textrm{for} & y \in (x_2,\infty), \\
  \end{array}
\end{equation}
where the points $x_1 < x_2$ are determined by the condition $\frac{\partial s}{\partial y} = \frac{\partial i}{\partial y}$ for a fixed value of $r$; see \cite{barakova}.

\subsection{The Model via sign-based representation}
\label{sign-based_representation}

The sign-based model representation, used in trend-based analysis (see \cite{doubravsky_dohnal}, \cite{forbus}, \cite{vicha_dohnal}), is characterised by a set of signs with three possible states: positive~$+$, negative~$-$, and zero~$0$, assigned to each model element. This set of signs includes the sign of a variable or function value, the sign of first- and second-order time derivatives of model variables, and the sign of first- and second-order partial derivatives of model functions with respect to each variable. Specifically, a model variable $x$ is represented by a triplet ($X$~D$X$~DD$X$), where $X$ denotes the sign of the variable $x$, D$X$ denotes the sign of its first-order time derivative, and DD$X$ denotes the sign of its second-order time derivative. A model function $f$ is represented, with respect to each variable $x$, by a triplet ($F_x$~D$F_x$~DD$F_x$), where $F_x$ denotes the sign of the function value, D$F_x$ the sign of the first-order partial derivative of $f$, and DD$F_x$ the sign of the second-order partial derivative of $f$. A similar situation holds for a constant (understood as a constant function), with the difference that the first- and second-order partial derivatives are zero. Thus, a constant is represented by a triplet ($C~0~0$), where $C$ denotes the sign of the constant $c$. It is possible to work with an incomplete set of model signs in the sense that some of these signs may be unavailable.

In Table~\ref{tab:model_variables_constant_representations}, we present the sign-based representations of the \mbox{IS--LM} model variables $y>0$ and $r \in \mathbb{R}$, and the constants $m>0$ and $i_a - s_a > 0$. The symbol X denotes that all three signs $+,-,0$ are possible.
\begin{table}[ht]
\begin{tabular}{lcccccccc}
\hline
System                                  & $y$   & $r$   &  $m$      & $i_a-s_a$ \\
\hline
Linear model \eqref{linear_IS-LM_model} & $+$XX & $+$XX & $+\,0\,0$ & $+\,0\,0$ \\
Nonlinear model \eqref{IS-LM_model}     & $+$XX & XXX   & $+\,0\,0$ & \\
\hline \hline
\end{tabular}
\vspace{0.2cm}
\caption{Sign-based representations of the \mbox{IS--LM} model variables and constants}
\label{tab:model_variables_constant_representations}
\end{table}

Subsequently, such a system of ordinary differential equations can be rewritten using sign-based representations as follows. Each model element is replaced by its corresponding sign-based symbol, i.e., a model variable $x$ by $X$, a model function $f$ by $F$, and a constant $c$ by $C$, while the first-order time derivatives of $x$ are represented by D$X$. Positive and negative multiplicative parameters are represented by the signs $+$ and $-$, respectively. Moreover, the addition rules for sign representations are given in Table~\ref{tab:sign-based_addition}. We can see that the addition of positive and negative signs cannot be determined, i.e., it is unknown. Therefore, in sign-based representations of models described by systems of differential equations, such ambiguous additions are avoided by transferring terms with negative signs to the other side of the equation, so that only additions of signs remain in the representation.
\begin{table}[ht]
\begin{tabular}{c|ccc}
\hline
& $+$ & $0$ & $-$  \\
\hline
  $+$ & $+$ & $+$ & ? \\
  $0$ & $+$ & $0$ & $-$ \\
  $-$ &  ?  & $-$ & $-$ \\
\hline \hline
\end{tabular}
\vspace{0.2cm}
\caption{Addition rules for sign representations}
\label{tab:sign-based_addition}
\end{table}

In Table~\ref{tab:systems_representations}, the considered \mbox{IS--LM} models are expressed in their sign-based representations. We first present the linear system \eqref{linear_IS-LM_model}, denoted by L1, and three models derived from it, denoted by L2, L3, and L4, where $n := \sigma - a > 0$ in the case $a < \sigma$, and $p := a - \sigma > 0$ in the case $a > \sigma$. Similarly, we present the \mbox{IS--LM} model in its general form \eqref{IS-LM_model}, denoted by N1, together with a simplified version, denoted by N2, corresponding to a particular subclass of models where $i(y,r) = i(y) + i(r)$, $s(y,r) = s(y) + s(r)$, and $i(r) - s(r) := q(r)$. Thus, $DY$ and $DR$ denote $\frac{dy}{dt}$ and $\frac{dr}{dt}$, $Y$ and $R$ denote $y$ and $r$, and $I, S, L, Q, M,$ and $A$ denote $i, s, l, q, m,$ and $i_a - s_a$, respectively.
\begin{table}[!htbp]
\begin{tabular}{lll}
\hline
Notation & System of differential equations & Sign-based representation \\
\hline
L1 & $\frac{dy}{dt} = \alpha \bigl[ ay - \sigma y - (b+d) r + (i_a - s_a) \bigr]$ & $DY + Y + R = Y + A$ \\
& $\frac{dr}{dt} = \beta \bigl[ k y - h r - m \bigr]$ & $DR + R + M = Y$ \\
\hline
L2 & $\frac{dy}{dt} = \alpha \bigl[ -n y - (b+d) r + (i_a - s_a) \bigr]$ & $DY + Y + R = A$ \\
& $\frac{dr}{dt} = \beta \bigl[ k y - h r - m \bigr]$ & $DR + R + M = Y$ \\
\hline
L3 & $\frac{dy}{dt} = \alpha \bigl[ -(b+d) r + (i_a - s_a) \bigr]$ & $DY + R = A$ \\
& $\frac{dr}{dt} = \beta \bigl[ k y - h r - m \bigr]$ & $DR + R + M = Y$ \\
\hline
L4 & $\frac{dy}{dt} = \alpha \bigl[ p y - (b+d) r + (i_a - s_a) \bigr]$ & $DY + R = Y + A$ \\
& $\frac{dr}{dt} = \beta \bigl[ k y - h r - m \bigr]$ & $DR + R + M = Y$ \\
\hline
N1 & $\frac{dy}{dt} = \alpha \bigl[ i(y,r) - s(y,r) \bigr]$ & $DY + S = I$ \\
& $\frac{dr}{dt} = \beta \bigl[ l(y,r) - m \bigr]$ & $DR + M = L$ \\
\hline
N2 & $\frac{dy}{dt} = \alpha \bigl[ i(y) - s(y) + q(r) \bigr]$ & $DY + S = I + Q$ \\
& $\frac{dr}{dt} = \beta \bigl[ l(y,r) - m \bigr]$ & $DR + M = L$ \\
\hline \hline 
\end{tabular}
\vspace{0.2cm}
\caption{Sign-based representations of the \mbox{IS--LM} model differential equations}
\label{tab:systems_representations}
\end{table}

Thus, the sign-based representation of the linear models consists of the variable and constant representations listed in Table~\ref{tab:model_variables_constant_representations} corresponding to the linear model, together with the system representations L1, L2, L3, and L4 given in Table~\ref{tab:systems_representations}.

\begin{figure}[ht]
  \centering
  \includegraphics[height=3.2cm]{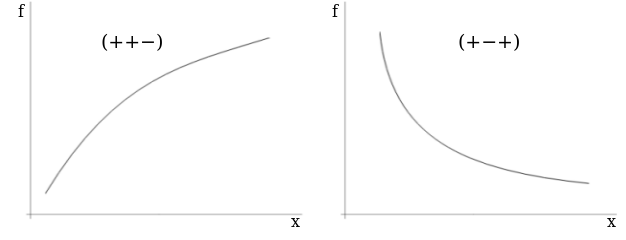}
  \caption{Illustration of decelerating increase and decelerating decrease}
  \label{fig:decelerating_inc_dec}
\end{figure}

For nonlinear models, a sign-based representation of the properties of the model functions is also required. Thus, we describe these properties using triplets ($F_x$~D$F_x$~DD$F_x$). The economic properties \eqref{economic_is_y}, \eqref{economic_is_r}, and \eqref{economic_l_yr} provide information only about $I_y$, D$I_y$, $S_y$, D$S_y$, $I_r$, D$I_r$, $S_r$, D$S_r$, $L_y$, D$L_y$, and $L_r$, D$L_r$. To obtain information about the second-order derivatives DD$I_r$, DD$S_r$, DD$L_y$, and DD$L_r$, we additionally assume decelerating increase or decrease, in the sense illustrated in Figure~\ref{fig:decelerating_inc_dec}.

Thus, the sign-based representation of the property \eqref{economic_is_r} with decelerating decrease for the investment function is ($I_r$~D$I_r$~DD$I_r$) = $(+-+)$, and with decelerating increase for the saving function is ($S_r$~D$S_r$~DD$S_r$) = $(++-)$. The properties of the money demand function \eqref{economic_l_yr} with decelerating increase and decrease are expressed by ($L_y$~D$L_y$~DD$L_y$) = $(++-)$ and ($L_r$~D$L_r$~DD$L_r$) = $(+-+)$, respectively. In addition, the Kaldor conditions \eqref{Kaldor_is}, in conjunction with \eqref{economic_is_y}, describe an accelerating increase of the investment and saving functions with respect to $y$ over certain intervals, namely up to $y = z_1$ for the investment function and from $y = z_2$ for the saving function, see Figure~\ref{fig:Kaldor_is}. Thus, the Kaldor conditions can be reformulated in the sense of sign-based representation as given in \eqref{Kaldor_i} and \eqref{Kaldor_s}.
\begin{equation}
\label{Kaldor_i}
  \begin{array}{lll}
  (I_y, \text{D}I_y, \text{DD}I_y)=(+++) & \textrm{for} & y \in (-\infty,z_1), \\
  (I_y, \text{D}I_y, \text{DD}I_y)=(++-) & \textrm{for} & y \in (z_1,\infty),
  \end{array}
\end{equation}
\begin{equation}
\label{Kaldor_s} 
  \begin{array}{lll} 
  (S_y, \text{D}S_y, \text{DD}S_y)=(++-) & \textrm{for} & y \in (-\infty,z_2), \\
  (S_y, \text{D}S_y, \text{DD}S_y)=(+++) & \textrm{for} & y \in (z_2,\infty), \\
  \end{array}
\end{equation}
where $z_1, z_2 \in (x_1,x_2)$.

In Table~\ref{tab:model_functions_representations}, we summarize the sign-based representations of the model functions with respect to each variable, namely $i(y)$, $s(y)$, $i(r)$, $q(r)$, $l(y)$, and $l(r)$. This table also provides the list of nonlinear models to which the trend-based analysis is applied. The first column contains the model notation, and the second column contains the system used in the model. The same system can be used in different models, which differ in their model function properties represented by the corresponding sign triplets. The system N1 is divided into three models. Model N1.1 represents the case of monotonically increasing functions $i(y)$ and $s(y)$. Models N1.2 and N1.3 incorporate the Kaldor conditions, distinguishing between two branches: the first with accelerating increase of $i(y)$ over the interval $(0,z_1)$, and the second with accelerating increase of $s(y)$ over the interval $(z_2,\infty)$. In these models, the property ($S_r$~D$S_r$~DD$S_r$) is omitted due to computational considerations. Therefore, we introduce the models N2.1, N2.2, and N2.3, in which this property is incorporated through the function $q(r)=i(r)-s(r)$. The sign-based representation of $q(r)$ follows from the subtraction of these functions, with three possible sign configurations of its value arising from this operation, see Table~\ref{tab:sign-based_addition}.
\begin{table}[ht]
\begin{tabular}{llcccccc}
\hline
Model & System &  $i(y)$    & $s(y)$  & $i(r)$ & $q(r)$ & $l(y)$ & $l(r)$ \\
\hline
N1.1  & N1     & $+\,+$\,X  & $+\,+$\,X  & $+-+$  &         & $++-$ & $+-+$ \\
N1.2  & N1     & $+++$      & $++-$      & $+-+$  &         & $++-$ & $+-+$ \\
N1.3  & N1     & $++-$      & $+++$      & $+-+$  &         & $++-$ & $+-+$ \\
N2.1  & N2     & $+\,+$\,X  & $+\,+$\,X  &        & $+-+$   & $++-$ & $+-+$ \\
N2.2  & N2     & $+\,+$\,X  & $+\,+$\,X  &        & $0\,-+$ & $++-$ & $+-+$ \\
N2.3  & N2     & $+\,+$\,X  & $+\,+$\,X  &        & $--+$   & $++-$ & $+-+$ \\
\hline \hline
\end{tabular}
\vspace{0.2cm}
\caption{Sign-based representations of the \mbox{IS--LM} model function properties}
\label{tab:model_functions_representations}
\end{table}

Thus, the sign-based representation of the nonlinear models N1.1, N1.2, N1.3, N2.1, N2.2, and N2.3 consists of the variable and constant representations listed in Table~\ref{tab:model_variables_constant_representations} corresponding to the nonlinear model, together with the system representations N1 and N2 listed in Table~\ref{tab:systems_representations} and the function representations listed in Table~\ref{tab:model_functions_representations}.

\section{Qualitative Analysis of the Model}

\subsection{Phase-Plane Qualitative Analysis of the Model}
\label{phase_plane_analysis}

The dynamics of a two-dimensional \mbox{IS--LM} model have been extensively studied in the general nonlinear form, both in Kaldor-type and non-Kaldor-type systems; see, e.g., \cite{chang_smyth}, \cite{gandolfo}, \cite{shone}, \cite{torre}, \cite{varian}, \cite{volna}.

The phase-plane qualitative analysis, based on linearization around stationary points and the classification of singular points, employs the Jacobian matrix, i.e., the matrix of first-order partial derivatives of the right-hand side functions of the system. Thus, the model dynamics can be inferred from properties \eqref{economic_is_y}, \eqref{economic_is_r}, \eqref{economic_l_yr}, and \eqref{Kaldor_is}. Therefore, we provide a brief overview and analysis of model singularities for the linear model \eqref{linear_IS-LM_model} and the nonlinear model \eqref{IS-LM_model} under economic conditions \eqref{economic_is_y}, \eqref{economic_is_r}, and \eqref{economic_l_yr}, with the nonlinear Kaldor-type system additionally satisfying \eqref{Kaldor_is}, serving as a representative illustration of model dynamics. The considered cases capture the essential qualitative behaviour of the system within the phase-plane framework, with other possible configurations representing variations of the same underlying dynamics.

In the \mbox{IS--LM} model, the nullclines are referred to as the IS curve and the LM curve. The IS and LM curves consist of points corresponding to equilibrium in the goods and money markets, respectively. According to the Implicit Function Theorem, the slopes of the IS and LM curves are given by
\begin{equation}
\label{IS_LM_curves_slopes}
\left.\frac{dr}{dy}\right|_{IS}=-\frac{i_y-s_y}{i_r-s_r},\quad \left.\frac{dr}{dy}\right|_{LM}=-\frac{l_y}{l_r},
\end{equation}
at points on the IS and LM curves, respectively, where $i_y, i_r, s_y, s_r, l_y$ and $l_r$ denote the corresponding partial derivatives. For the linear case \eqref{linear_IS-LM_model}, the equations of the IS and LM curves are given by
\begin{equation}
\label{IS_LM_curves_linear}
\begin{array}{llcl}
\text{IS curve:} & r &=& \frac{a-\sigma}{b+d}y + \frac{i_a-s_a}{b+d}, \\
\text{LM curve:} & r &=& \frac{k}{h}y + \frac{m}{h} .
\end{array}
\end{equation}

Each intersection point of the IS and LM curves corresponds to a stationary point, representing aggregate macroeconomic equilibrium, i.e., the simultaneous equilibrium in the goods and money markets. The analysis is restricted to cases where the IS and LM curves intersect, as non-intersecting configurations are of limited economic relevance. The linear model \eqref{linear_IS-LM_model} admits a single equilibrium point, except in the non-standard situation where the IS and LM curves overlap. Nonlinear models can exhibit more singular points, especially in the Kaldor-type system, where at least one and at most three singular points emerge, except in the non-standard situation in which the IS and LM curves overlap over some interval.

The core of the phase-plane qualitative analysis lies in the classification of singular points through the spectral analysis of the Jacobian matrix of the linearised system at the equilibrium point. Thus, the eigenvalues of the Jacobian matrix $J$ associated with the system \eqref{IS-LM_model} are given by
\begin{equation}
\label{eigenvalues}
\lambda_{1,2} = \frac{1}{2}\left[\operatorname{tr}J \pm \sqrt{(\operatorname{tr}J)^2 - 4\det J}\right],
\end{equation}
where $\operatorname{tr}J = \alpha (i_y - s_y) + \beta l_r$ and $\det J = \alpha \beta \left[ (i_y - s_y)l_r - (i_r - s_r)l_y \right]$. In particular, for the linear system \eqref{linear_IS-LM_model}, the trace and determinant of the corresponding Jacobian matrix are $\operatorname{tr}J = \alpha (a - \sigma) - \beta h$ and $\det J = \alpha \beta \left[ k(b+d)-h(a-\sigma) \right]$. Therefore, the phase-plane qualitative analysis is summarised in Tables \ref{tab:hyperbolic_eq_points}, \ref{tab:non-hyperbolic_eq_points}, \ref{tab:overlapping_curves_eq_points}, and \ref{tab:Kaldor_IS-LM_model_eq_points}. The entries in these tables are derived from \eqref{IS_LM_curves_slopes} and \eqref{eigenvalues}, based on the economic conditions \eqref{economic_is_y}, \eqref{economic_is_r}, and \eqref{economic_l_yr}, and in the Kaldor case also \eqref{Kaldor_is}, together with comparisons of the partial derivatives $i_y$ and $s_y$ and of the slopes of the IS and LM curves. From \eqref{economic_l_yr} and \eqref{IS_LM_curves_slopes}, it follows that the LM~curve is increasing over the entire considered interval for $y>0$. In accordance with the usual terminology, stable nodes and foci are called sinks, while unstable nodes and foci are called sources.

\begin{table}[!htbp]
\centering
\begin{tabular}{lll p{3.2cm} ccl}
\hline
Properties & IS curve & Slopes & Figure & $\det J$ & $\operatorname{tr}J$ & Type\\
\hline
$a<\sigma$ & decreasing & $\frac{a-\sigma}{b+d}<\frac{k}{h}$ &
  \multirow{2}{*}{\raisebox{-0.9\height}{\includegraphics[height=2.9cm]{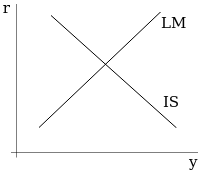}}}
  & $+$ & $-$ & sink \\
$i_y<s_y$ && $\frac{dr}{dy}\big|_{IS}<\frac{dr}{dy}\big|_{LM}$ \\[10.5ex]
\hline
$a=\sigma$ & constant & $\frac{a-\sigma}{b+d}=0<\frac{k}{h}$ &
  \multirow{2}{*}{\raisebox{-0.9\height}{\includegraphics[height=2.9cm]{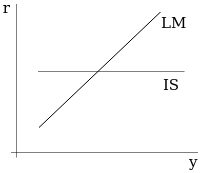}}}
  & $+$ & $-$ & sink \\
$i_y=s_y$ && $\frac{dr}{dy}\big|_{IS}=0<\frac{dr}{dy}\big|_{LM}$ \\[10.5ex]
\hline
$a>\sigma$ & increasing & $\frac{a-\sigma}{b+d}<\frac{k}{h}$ &
  \multirow{2}{*}{\raisebox{-0.9\height}{\includegraphics[height=2.9cm]{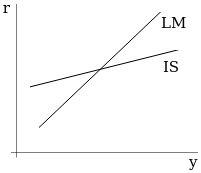}}}
  & $+$ & $-$ & sink \\
$i_y>s_y$ && $\frac{dr}{dy}\big|_{IS}<\frac{dr}{dy}\big|_{LM}$ && $+$ & $+$ & source \\[10.5ex]
\hline
$a>\sigma$ & increasing & $\frac{a-\sigma}{b+d}>\frac{k}{h}$ &
  \multirow{2}{*}{\raisebox{-0.9\height}{\includegraphics[height=2.9cm]{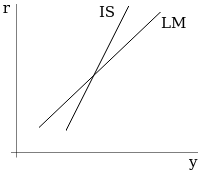}}}
  & $-$ & $-$ & saddle \\
$i_y>s_y$ && $\frac{dr}{dy}\big|_{IS}>\frac{dr}{dy}\big|_{LM}$ && $-$ & $+$ & saddle \\ 
          &&&& $-$ & $0$ & saddle \\[7ex]
\hline \hline
\end{tabular}
\vspace{0.1cm}
\caption{Classification of hyperbolic singular points in the linear and corresponding nonlinear \mbox{IS--LM} model under the given economic conditions in the case of a single equilibrium}
\label{tab:hyperbolic_eq_points}
\end{table}

In Table~\ref{tab:hyperbolic_eq_points}, we summarize the classification of hyperbolic singular points in the linear \mbox{IS--LM} model \eqref{linear_IS-LM_model}. This classification can also be applied to the nonlinear \mbox{IS--LM} model \eqref{IS-LM_model} under economic conditions \eqref{economic_is_y}, \eqref{economic_is_r}, and \eqref{economic_l_yr} in the case of a single hyperbolic equilibrium, regarded as a classification of the corresponding linearised system at this equilibrium point. The figures included in Table~\ref{tab:hyperbolic_eq_points} illustrate linear cases or the considered nonlinear cases in terms of the corresponding linearised systems in a neighbourhood of the singular point.

\begin{table}[ht]
\centering
\begin{tabular}{ll p{6.5cm} ccl}
\hline
System & Slopes & Figure & $\det J$ & $\operatorname{tr}J$ & Type\\
\hline
Linear & $\frac{a-\sigma}{b+d}<\frac{k}{h}$ &
  \multirow{2}{*}{\raisebox{-0.9\height}{\includegraphics[height=2.9cm]{figures/linear_IS_LM_curves_3.png}}}
  & $+$ & $0$ & centre \\[14ex]
\hline
Nonlinear & $\frac{dr}{dy}\big|_{IS}<\frac{dr}{dy}\big|_{LM}$ &  \multirow{2}{*}{\raisebox{-0.9\height}{\includegraphics[height=2.9cm]{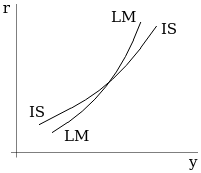}}} & $+$ & $0$ & N/A \\[14ex]
\hline
Nonlinear & $\frac{dr}{dy}\big|_{IS}=\frac{dr}{dy}\big|_{LM}$ &  \multirow{2}{*}{\raisebox{-0.9\height}{\includegraphics[height=2.9cm]{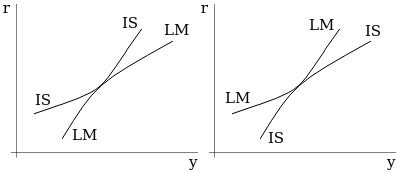}}} & $0$ & $-$ & N/A \\
  &&&& $+$ & N/A \\ 
  &&&& $0$ & N/A \\[7ex]
\hline \hline
\end{tabular}
\vspace{0.1cm}
\caption{Classification of non-hyperbolic singular points in the linear and corresponding nonlinear \mbox{IS--LM} model under the given economic conditions in the case of a single equilibrium}
\label{tab:non-hyperbolic_eq_points}
\end{table}

In Table~\ref{tab:non-hyperbolic_eq_points}, we illustrate the analysis of \mbox{IS--LM} systems characterised by a single non-hyperbolic singular point, such as a centre or degenerate cases of tangent IS and LM curves. In the considered setting, non-hyperbolic singular points may arise in the case where $i_y > s_y$, which in the linear system corresponds to $a>\sigma$. This situation leads to an increasing IS curve. For non-hyperbolic singular points, analysis based on linearisation and spectral analysis of the Jacobian matrix does not provide sufficient information about the system dynamics. The figures included in Table~\ref{tab:non-hyperbolic_eq_points} illustrate cases with a single non-hyperbolic singular point in the linear system and in the corresponding nonlinear system.

\begin{table}[ht]
\centering
\begin{tabular}{ll p{3.2cm} cl}
\hline
System & Slopes & Figure & $\operatorname{tr}J$ & Type\\
\hline
Linear & $\frac{a-\sigma}{b+d}=\frac{k}{h}$ &
  \multirow{2}{*}{\raisebox{-0.9\height}{\includegraphics[height=2.9cm]{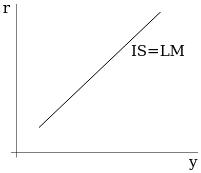}}}
  & $-$ & attracting set \\ 
&&& $+$ & repelling set \\ 
&&& $0$ & neutral (degenerate) \\[7ex]
\hline
Nonlinear & $\frac{dr}{dy}\big|_{IS}=\frac{dr}{dy}\big|_{LM}$ &
  \multirow{2}{*}{\raisebox{-0.9\height}{\includegraphics[height=2.9cm]{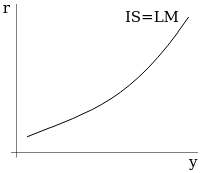}}}
  & $-$ & N/A \\ 
&&& $+$ & N/A \\ 
&&& $0$ & N/A \\[7ex]
\hline \hline
\end{tabular}
\vspace{0.1cm}
\caption{Singular points in the linear and corresponding nonlinear \mbox{IS--LM} model under the given economic conditions in the case of overlapping IS and LM curves}
\label{tab:overlapping_curves_eq_points}
\end{table}

In Table~\ref{tab:overlapping_curves_eq_points}, we provide a qualitative description of \mbox{IS--LM} systems characterised by infinitely many non-hyperbolic singular points arising when the IS and LM curves, both increasing, overlap in the linear \mbox{IS--LM} model \eqref{linear_IS-LM_model} and in the corresponding nonlinear system~\eqref{IS-LM_model}. In the considered setting, such configurations arise when $i_y > s_y$, especially for the linear system with $a>\sigma$, and when the slopes of the IS and LM curves are equal over the entire considered interval in the linear case and over at least part of the considered interval in the nonlinear case. Equal slopes of the IS and LM curves over this interval imply a zero determinant of the corresponding Jacobian matrix. In this case, a set of equilibria forms a line described by~\eqref{IS_LM_curves_linear} in the linear case, and a curve in the nonlinear case. In such non-standard cases, analysis based on linearisation and spectral analysis of the Jacobian matrix does not provide sufficient information about the system dynamics. The figures included in Table~\ref{tab:overlapping_curves_eq_points} illustrate overlapping IS and LM curves for the linear system and the corresponding nonlinear system.

\begin{table}[!htbp]
\centering
\begin{tabular}{p{5cm} ccc}
\hline
Global configuration & $y \in (0,x_1]$ & $y \in (x_1,x_2)$ & $y \in [x_2, \infty)$ \\
\hline
\multirow{2}{*}{\raisebox{-0.8\height}{\includegraphics[height=1.5cm]{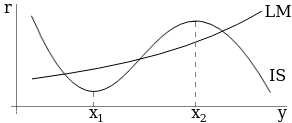}}} &
\multirow{2}{*}{\raisebox{-0.5\height}{\includegraphics[height=1cm]{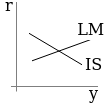}}} &
\multirow{2}{*}{\raisebox{-0.5\height}{\includegraphics[height=1cm]{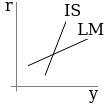}}} &
\multirow{2}{*}{\raisebox{-0.5\height}{\includegraphics[height=1cm]{figures/Kaldor_sink_1.png}}} \\[2.4ex]
& sink & saddle & sink \\
\hline
\multirow{2}{*}{\raisebox{-0.8\height}{\includegraphics[height=1.5cm]{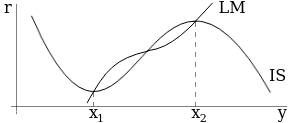}}} &
\multirow{2}{*}{\raisebox{-0.5\height}{\includegraphics[height=1cm]{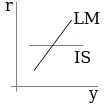}}} &
\multirow{2}{*}{\raisebox{-0.5\height}{\includegraphics[height=1cm]{figures/Kaldor_saddle.png}}} &
\multirow{2}{*}{\raisebox{-0.5\height}{\includegraphics[height=1cm]{figures/Kaldor_sink_2.png}}} \\[2.4ex] 
& sink & saddle & sink \\ 
\hline
\multirow{2}{*}{\raisebox{-0.8\height}{\includegraphics[height=1.5cm]{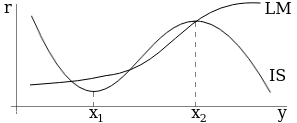}}} &
\multirow{2}{*}{\raisebox{-0.5\height}{\includegraphics[height=1cm]{figures/Kaldor_sink_1.png}}} &
\multirow{2}{*}{\raisebox{-0.5\height}{\includegraphics[height=1cm]{figures/Kaldor_saddle.png}}} &
\multirow{2}{*}{\raisebox{-0.5\height}{\includegraphics[height=1cm]{figures/Kaldor_sink_2.png}}} \\[2.4ex]
& sink & saddle & sink \\ 
\hline
\multirow{2}{*}{\raisebox{-0.8\height}{\includegraphics[height=1.5cm]{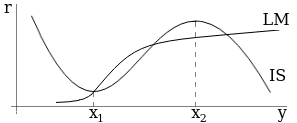}}} &
\multirow{2}{*}{\raisebox{-0.5\height}{\includegraphics[height=1cm]{figures/Kaldor_sink_2.png}}} &
\multirow{2}{*}{\raisebox{-0.5\height}{\includegraphics[height=1cm]{figures/Kaldor_saddle.png}}} &
\multirow{2}{*}{\raisebox{-0.5\height}{\includegraphics[height=1cm]{figures/Kaldor_sink_1.png}}} \\[2.4ex]
& sink & saddle & sink \\ 
\hline
\multirow{2}{*}{\raisebox{-0.8\height}{\includegraphics[height=1.5cm]{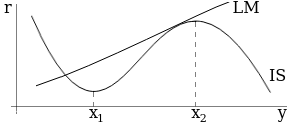}}} &
\multirow{2}{*}{\raisebox{-0.5\height}{\includegraphics[height=1cm]{figures/Kaldor_sink_1.png}}} &
\multirow{2}{*}{\raisebox{-0.5\height}{\includegraphics[height=1cm]{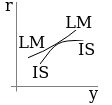}}} \\[2.4ex]
& sink & N/A & \\ 
\hline
\multirow{2}{*}{\raisebox{-0.8\height}{\includegraphics[height=1.5cm]{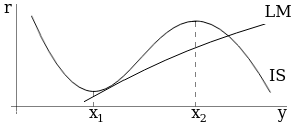}}} & &
\multirow{2}{*}{\raisebox{-0.5\height}{\includegraphics[height=1cm]{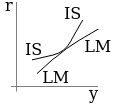}}} &
\multirow{2}{*}{\raisebox{-0.5\height}{\includegraphics[height=1cm]{figures/Kaldor_sink_1.png}}} \\[2.4ex]
&  & N/A & sink \\ 
\hline
\multirow{2}{*}{\raisebox{-0.8\height}{\includegraphics[height=1.5cm]{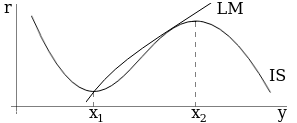}}} &
\multirow{2}{*}{\raisebox{-0.5\height}{\includegraphics[height=1cm]{figures/Kaldor_sink_2.png}}} &
\multirow{2}{*}{\raisebox{-0.5\height}{\includegraphics[height=1cm]{figures/Kaldor_tangent_curves_1.png}}} \\[2.4ex]
& sink & N/A &  \\ 
\hline
\multirow{2}{*}{\raisebox{-0.8\height}{\includegraphics[height=1.5cm]{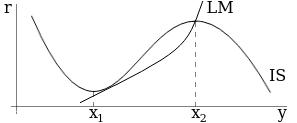}}} & &
\multirow{2}{*}{\raisebox{-0.5\height}{\includegraphics[height=1cm]{figures/Kaldor_tangent_curves_2.png}}} &
\multirow{2}{*}{\raisebox{-0.5\height}{\includegraphics[height=1cm]{figures/Kaldor_sink_2.png}}} \\[2.4ex]
&  & N/A & sink \\ 
\hline
\multirow{2}{*}{\raisebox{-0.8\height}{\includegraphics[height=1.5cm]{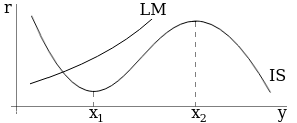}}} &
\multirow{2}{*}{\raisebox{-0.5\height}{\includegraphics[height=1cm]{figures/Kaldor_sink_1.png}}} \\[2.4ex]
& sink &  &  \\ 
\hline
\multirow{2}{*}{\raisebox{-0.8\height}{\includegraphics[height=1.5cm]{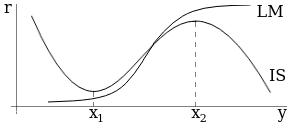}}} & &
\multirow{2}{*}{\raisebox{-0.5\height}{\includegraphics[height=1cm]{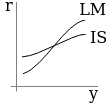}}} \\[2.4ex]
& & sink, source, N/A \\
\hline
\multirow{2}{*}{\raisebox{-0.8\height}{\includegraphics[height=1.5cm]{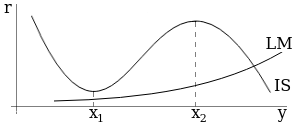}}} & & &
\multirow{2}{*}{\raisebox{-0.5\height}{\includegraphics[height=1cm]{figures/Kaldor_sink_1.png}}} \\[2.4ex]
& & & sink \\
\hline
\multirow{2}{*}{\raisebox{-0.8\height}{\includegraphics[height=1.5cm]{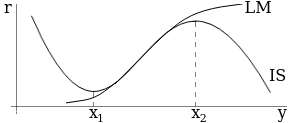}}} & &
\multirow{2}{*}{\raisebox{-0.5\height}{\includegraphics[height=1cm]{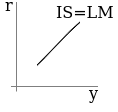}}} & \\[2.4ex]
& & N/A & \\ 
\hline \hline
\end{tabular}
\vspace{0.1cm}
\caption{Classification of singular points in the Kaldor-type \mbox{IS--LM} model}
\label{tab:Kaldor_IS-LM_model_eq_points}
\end{table}

Finally, in Table~\ref{tab:Kaldor_IS-LM_model_eq_points}, we summarize the phase-plane qualitative analysis for the Kaldor-type \mbox{IS--LM} model \eqref{IS-LM_model} under the conditions \eqref{economic_is_y}, \eqref{economic_is_r}, \eqref{economic_l_yr}, and \eqref{Kaldor_is}. For this classification, we build on the cases described in Tables \ref{tab:hyperbolic_eq_points}, \ref{tab:non-hyperbolic_eq_points}, and \ref{tab:overlapping_curves_eq_points}. The first column presents global configurations of the IS and LM curves, from which the number and positions of singular points can be identified. In the Kaldor-type \mbox{IS--LM} model, the IS curve is decreasing for $y \in (0,x_1) \cup (x_2, \infty)$ and increasing for $y \in (x_1,x_2)$. The other columns illustrate the local situation through the corresponding linearised system in the case of hyperbolic singular points, and a zoomed local representation in the case of non-hyperbolic singular points or overlapping IS and LM curves, corresponding to the appropriate interval. In these columns, under each figure illustrating the local situation, the type of hyperbolic singular points is identified, while non-hyperbolic cases are indicated by “N/A”, as shown in Tables \ref{tab:hyperbolic_eq_points}, \ref{tab:non-hyperbolic_eq_points}, and \ref{tab:overlapping_curves_eq_points}.

As can be seen, phase-plane qualitative analysis based on the classification of singular points and limited knowledge of first-order partial derivatives of the model functions provides a complete description in the case of hyperbolic singular points. For other cases, this method does not allow for a conclusive classification and additional information about the model functions, such as their types (e.g. polynomial, trigonometric, exponential), is required in order to describe the system dynamics. Thus, the classification of non-hyperbolic singular points depends strongly on additional information about higher-order derivatives of the model functions.

\subsection{Trend-Based Qualitative Analysis of the Model}
\label{trend_based_qualitative_analysis}

A detailed description of the principles of trend-based analysis can be found in \cite{dohnal_doubravsky_15} and \cite{dohnal_doubravsky_16}. In general, trend-based analysis follows from the sign-based representation specified in Section~\ref{sign-based_representation} and generates a set of possible scenarios together with possible transitions between them for the considered models.

A scenario is defined by the sign triplets assigned to each model variable. Formally, the set of $m$ possible scenarios corresponding to $n$ model variables $X_i$ is given by
\begin{equation}
S(n, m) = \{\{(X_1~\text{D}X_1~\text{DD}X_1), (X_2~\text{D}X_2~\text{DD}X_2 ), \dots, (X_n~\text{D}X_n~\text{DD}X_n)\}_j \},
\end{equation}
where $i=1,\dots,n$ and $j=1,\dots,m$. The scenario $\{(X_1~0~0), (X_2~0~0), \dots, (X_n~0~¨0)\}$ represents a steady state of the model, i.e., a possible equilibrium configuration. Besides the set of possible scenarios corresponding to a model, a set of possible transitions between them, denoted by $T$, is also required. The identification of possible transitions is based on mathematical rules associated with model variables, which are typically assumed to be continuous and smooth time functions in autonomous systems of differential equations.

The set of all possible scenarios $S$ represents all possible configurations of model variables in terms of positivity, monotonicity, and curvature of the corresponding time functions. The set of all possible transitions $T$ represents all possible time evolutions of the model variables, both in the future and in the past.

The graphical representation of these results is given by a transition graph
\begin{equation}
H = (S, T)
\end{equation}
where the scenario set $S$ and the transition set $T$ form the nodes and directed edges of the directed graph $H$, respectively.

To generate the scenario and transition sets corresponding to the considered models specified in Section~\ref{sign-based_representation}, we use a combinatorial algorithm based on an artificial intelligence approach, for details, see \cite{dohnal_doubravsky_15}, \cite{dohnal_doubravsky_16}. This algorithm follows from the sign-based representation of the models specified in Tables~\ref{tab:model_variables_constant_representations}, \ref{tab:systems_representations},~and~\ref{tab:model_functions_representations}.

In the following, for each considered model, we present a table of scenarios and transition graphs together with accompanying explanations. In the table headings, the symbols $Y$ and $R$ represent the sign triplets ($Y$~D$Y$~DD$Y$) and ($R$~D$R$~DD$R$), respectively.

\subsubsection{Linear Models}
The linear \mbox{IS--LM} models L1, L2, L3, and L4 are considered in their sign-based qualitative representation. Each model shares the following characteristics in the sense of symmetric sets of scenarios and transitions. The first position of each triplet corresponding to $Y$ and $R$ is fixed as $+$ due to the positivity condition, see Table~\ref{tab:model_variables_constant_representations}. Under the adopted scenario numbering, pairs of scenarios no.~$s \pm k$, where $s$ denotes the number assigned to the steady-state scenario and $k = 1,2,\dots,s-1$, differ only by opposite signs in the second and third positions of the corresponding triplets. This symmetry also occurs in the transition sets and is apparent from the transition graphs, with the steady-state scenario located in the middle.

\begin{figure}[ht]
  \centering
  \includegraphics[height=7.7cm]{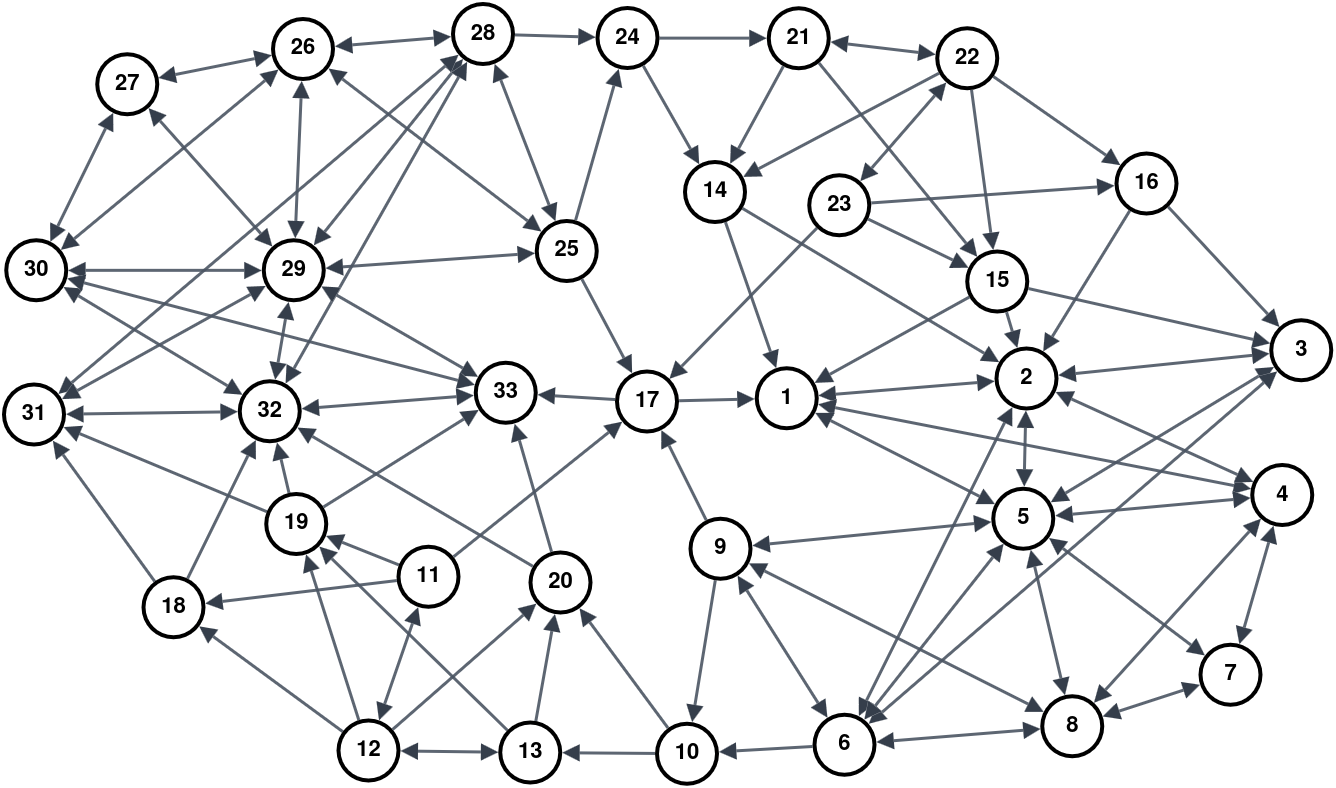}
  \caption{Transition graph corresponding to L1}
  \label{fig:transition_graph_L1}
\end{figure}

The scenario set corresponding to L1, listed in Table~\ref{tab:scenarios_L1}, consists of 33 scenarios and contains the steady-state scenario no.~17. The transition graph of L1 is shown in Figure~\ref{fig:transition_graph_L1}. However, this transition graph is not sufficiently transparent for observing the admissible properties of the system. Therefore, we decompose the model L1 into the three models L2, L3, and L4, see Table~\ref{tab:systems_representations}. The transition graphs of these models can be regarded as decomposed parts of the transition graph of L1, preserving the qualitative information contained in L1 in a more transparent form.

\begin{figure}[!htbp]
\centering
\begin{minipage}{0.48\textwidth}
    \centering
    \includegraphics[width=\textwidth]{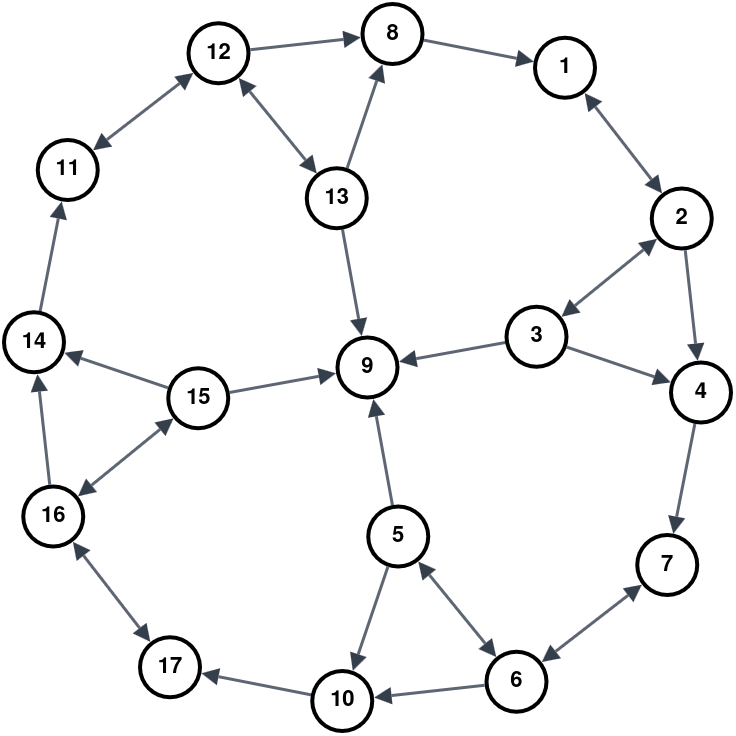}
    \caption{Transition graph corresponding to~L2}
    \label{fig:transition_graph_L2}
\end{minipage}
\hfill
\begin{minipage}{0.48\textwidth}
    \centering
    \includegraphics[width=\textwidth]{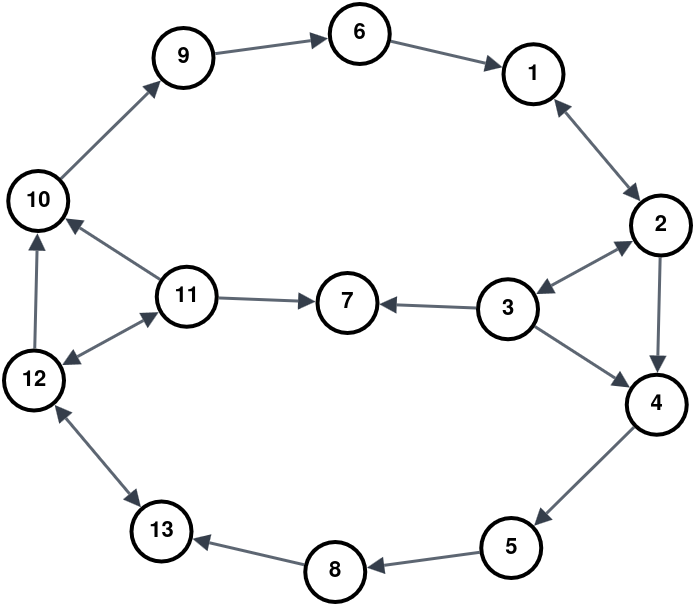}
    \caption{Transition graph corresponding to~L3}
    \label{fig:transition_graph_L3}
\end{minipage}
\end{figure}

\begin{figure}[!htbp]
  \centering
  \includegraphics[height=7.7cm]{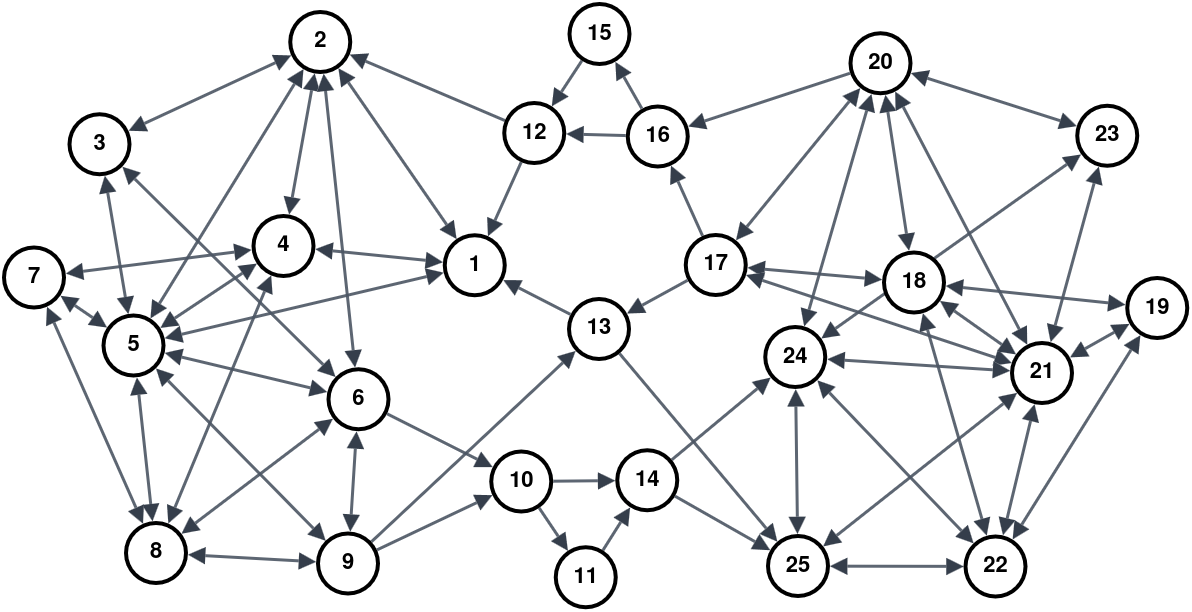}
  \caption{Transition graph corresponding to L4}
  \label{fig:transition_graph_L4}
\end{figure}

The scenario sets corresponding to L2, L3, and L4 consist of 17, 13, and 25 scenarios and contain steady-state scenarios no.~9, 7, and 13, respectively. These scenario sets are listed in Tables~\ref{tab:scenarios_L2}, \ref{tab:scenarios_L3}, and \ref{tab:scenarios_L4}. In the transition graphs corresponding to L2, L3, and L4, shown in Figures~\ref{fig:transition_graph_L2}, \ref{fig:transition_graph_L3}, and \ref{fig:transition_graph_L4}, respectively, we can identify cycles. Examples are $1 \rightarrow 2 \rightarrow 4 \rightarrow 7 \rightarrow 6 \rightarrow 10 \rightarrow 17 \rightarrow 16 \rightarrow 14 \rightarrow 11 \rightarrow 12 \rightarrow 8 \rightarrow 1$ for L2, $1 \rightarrow 2 \rightarrow 4 \rightarrow 5 \rightarrow 8 \rightarrow 13 \rightarrow 12 \rightarrow 10 \rightarrow 9 \rightarrow 6 \rightarrow 1$ for L3, and $13 \rightarrow 1 \rightarrow 4 \rightarrow 7 \rightarrow 8 \rightarrow 9 \rightarrow 13$ together with $13 \rightarrow 25 \rightarrow 22 \rightarrow 21 \rightarrow 23 \rightarrow 20 \rightarrow 17 \rightarrow 13$ for L4. The corresponding steady-state scenarios can be reached from scenarios no.~2, 6, 12, and 16 through scenarios no.~3, 5, 13, and 15 in the case of L2, and from scenarios no.~2 and 12 through scenarios no.~3 and 11 in the case of L3. In the case of L4, the steady-state scenario belongs to the identified cycles.

\begin{table}[!htbp]
\centering
\begin{tabular}{ccc c ccc c ccc}
\hline
No. & $Y$     & $R$     && No. & $Y$      & $R$      && No. & $Y$     & $R$ \\
\hline
1   & $+++$   & $+++$   && 12  & $+-\,0$  & $++-$    && 23  & $++-$   & $+-+$ \\
2   & $++\,0$ & $+++$   && 13  & $+--$    & $++-$    && 24  & $+~0~+$ & $+-+$ \\
3   & $++-$   & $+++$   && 14  & $+++$    & $+~0~+$  && 25  & $+-+$   & $+-+$ \\
4   & $+++$   & $++\,0$ && 15  & $++\,0$  & $+~0~+$  && 26  & $+-\,0$ & $+-+$ \\
5   & $++\,0$ & $++\,0$ && 16  & $++-$    & $+~0~+$  && 27  & $+--$   & $+-+$ \\
6   & $++-$   & $++\,0$ && 17  & $+~0~0~$ & $+~0~0~$ && 28  & $+-+$   & $+-\,0$ \\
7   & $+++$   & $++-$   && 18  & $+-+$    & $+~0~-$  && 29  & $+-\,0$ & $+-\,0$ \\
8   & $++\,0$ & $++-$   && 19  & $+-\,0$  & $+~0~-$  && 30  & $+--$   & $+-\,0$ \\
9   & $++-$   & $++-$   && 20  & $+--$    & $+~0~-$  && 31  & $+-+$   & $+--$ \\
10  & $+~0~-$ & $++-$   && 21  & $+++$    & $+-+$    && 32  & $+-\,0$ & $+--$ \\
11  & $+-+$   & $++-$   && 22  & $++\,0$  & $+-+$    && 33  & $+--$   & $+--$ \\
\hline \hline
\end{tabular}
\vspace{0.2cm}
\caption{Scenario set corresponding to L1}
\label{tab:scenarios_L1}
\end{table}

\begin{table}[ht]
\centering
\begin{tabular}{ccc c ccc c ccc}
\hline
No. & $Y$     & $R$     && No. & $Y$      & $R$      && No. & $Y$     & $R$ \\
\hline
1   & $++-$   & $+++$   && 7   & $+--$    & $++-$    && 12  & $++\,0$ & $+-+$ \\
2   & $++-$   & $++\,0$ && 8   & $++-$    & $+~0~+$  && 13  & $++-$   & $+-+$ \\
3   & $++-$   & $++-$   && \multirow{2}{*}{9} & \multirow{2}{*}{$+~0~0~$} & \multirow{2}{*}{$+~0~0~$}                                        && 14  & $+~0~+$ & $+-+$ \\
4   & $+~0~-$ & $++-$   &&     &&                    && 15  & $+-+$   & $+-+$ \\
5   & $+-+$   & $++-$   && 10  & $+-+$    & $+~0~-$  && 16  & $+-+$   & $+-\,0$ \\
6   & $+-\,0$ & $++-$   && 11  & $+++$    & $+-+$    && 17  & $+-+$   & $+--$ \\
\hline \hline
\end{tabular}
\vspace{0.2cm}
\caption{Scenario set corresponding to L2}
\label{tab:scenarios_L2}
\end{table}

\begin{table}[ht]
\centering
\begin{tabular}{ccc c ccc c ccc}
\hline
No. & $Y$     & $R$     && No. & $Y$      & $R$      && No. & $Y$     & $R$ \\
\hline
1   & $++-$   & $+++$   && 6   & $++\,0$  & $+~0~+$  && 9   & $+++$   & $+-+$ \\
2   & $++-$   & $++\,0$ &&     &          &          && 10  & $+~0~+$ & $+-+$ \\
3   & $++-$   & $++-$   && 7   & $+~0~0~$ & $+~0~0~$ && 11  & $+-+$   & $+-+$ \\
4   & $+~0~-$ & $++-$   &&     &          &          && 12  & $+-+$   & $+-\,0$\\
5   & $+--$   & $++-$   && 8   & $+-\,0$  & $+~0~-$  && 13  & $+-+$   & $+--$ \\
\hline \hline
\end{tabular}
\vspace{0.2cm}
\caption{Scenario set corresponding to L3}
\label{tab:scenarios_L3}
\end{table}

\newpage
\begin{table}[ht]
\centering
\begin{tabular}{ccc c ccc c ccc}
\hline
No. & $Y$     & $R$     && No. & $Y$      & $R$      && No. & $Y$     & $R$ \\
\hline
1 & $+++$   & $+++$   && 10 & $+~0~-$  & $++-$    && 17 & $+-+$   & $+-+$ \\
2 & $++\,0$ & $+++$   && 11 & $+--$    & $++-$    && 18 & $+-\,0$ & $+-+$\\
3 & $++-$   & $+++$   && 12 & $+++$    & $+~0~+$  && 19 & $+--$   & $+-+$ \\
4 & $+++$   & $++\,0$ &&    &          &          && 20 & $+-+$   & $+-\,0$\\
5 & $++\,0$ & $++\,0$ && 13 & $+~0~0~$ & $+~0~0~$ && 21 & $+-\,0$ & $+-\,0$ \\
6 & $++-$   & $++\,0$ &&    &          &          && 22 & $+--$   & $+-\,0$ \\
7 & $+++$   & $++-$   && 14 & $+--$    & $+~0~-$  && 23 & $+-+$   & $+--$ \\
8 & $++\,0$ & $++-$   && 15 & $+++$    & $+-+$    && 24 & $+-\,0$ & $+--$ \\ 
9 & $++-$   & $++-$   && 16 & $+~0~+$  & $+-+$    && 25 & $+--$   & $+--$ \\
\hline \hline
\end{tabular}
\vspace{0.2cm}
\caption{Scenario set corresponding to L4}
\label{tab:scenarios_L4}
\end{table}

\begin{rmk}
\label{rmk:linear_model}
From the perspective of trend-based modelling, the linear model can also be expressed using the sign-based representation of the general system N1 (see Table~\ref{tab:systems_representations}) together with the function property representations $(++0)$ and $(+-0)$ for increasing and decreasing linear functions, respectively. However, such a specification does not provide a sufficiently informative result. Due to insufficient structural information, applying the trend-based modelling algorithm to this specification results in only the steady-state scenario being admissible.
\end{rmk}

\subsubsection{Nonlinear Models}
The nonlinear models N1.1, N1.2, N1.3, N2.1, N2.2, and N2.3 are considered in their sign-based qualitative form. Although these models slightly differ in their system representations and model function properties, including the Kaldor conditions, see Tables~\ref{tab:systems_representations} and \ref{tab:model_functions_representations}, the trend-based modelling algorithm yields the same scenario and transition sets for all these nonlinear models. 

\begin{figure}[ht]
  \centering
  \includegraphics[height=6cm]{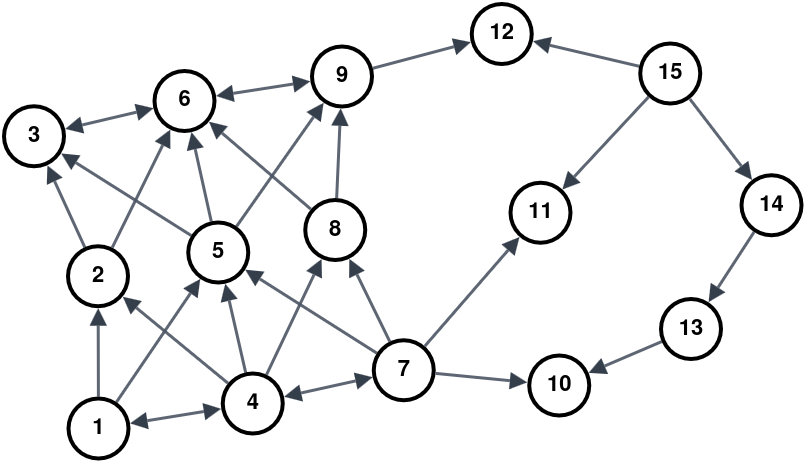}
  \caption{Transition graph corresponding to the nonlinear models}
  \label{fig:transition_graph_nonlinear}
\end{figure}

These scenario and transition sets do not exhibit the characteristic symmetry previously observed in the linear models. The first position of each triplet corresponding to $Y$ is fixed as $+$ due to the positivity condition, while the first position corresponding to $R$ varies among all three signs due to the absence of a positivity restriction on the interest rate in the nonlinear case, see Table~\ref{tab:model_variables_constant_representations}. The scenario set consists of 15 scenarios listed in Table~\ref{tab:scenarios_nonlinear} and contains three steady-state scenarios no.~10, 11, and 12, differing in the sign of the value of $R$. The corresponding transition graph is shown in Figure~\ref{fig:transition_graph_nonlinear}. From each node, there exists a directed path leading to one of the steady-state scenarios.
\begin{table}[ht]
\centering
\begin{tabular}{ccc c ccc c ccc}
\hline
No. & $Y$     & $R$     && No. & $Y$      & $R$       && No. & $Y$     & $R$ \\
\hline
1 & $+++$   & $+-+$     && 6   & $++\,0$  & $--+$     && 11 & $+~0~0$ & $\,0~0~0$ \\
2 & $+++$   & $\,0\,-+$ && 7   & $++-$    & $+-+$     && 12 & $+~0~0$ & $-\,0~0$ \\
3 & $+++$   & $--+$     && 8   & $++-$    & $\,0\,-+$ && 13 & $+-+$   & $++-$ \\ 
4 & $++\,0$ & $+-+$     && 9   & $++-$    & $--+$     && 14 & $+-+$   & $\,0\,+-$ \\ 
5 & $++\,0$ & $\,0\,-+$ && 10  & $+~0~0$  & $+~0~0$   && 15 & $+-+$   & $-+-$ \\ 
\hline \hline
\end{tabular}
\vspace{0.2cm}
\caption{Scenario set corresponding to the nonlinear models}
\label{tab:scenarios_nonlinear}
\end{table}

\begin{rmk}
As can be seen in Table~\ref{tab:model_functions_representations}, the nonlinear saving function property with respect to the interest rate, expressed in the sign-based representation as ($S_r$~D$S_r$~DD$S_r$) = $(++-)$, is not explicitly included in the considered nonlinear models. The trend-based modelling algorithm cannot consistently incorporate this property into the sign-based equation $DY + S = I$ due to sign incompatibility. The incompatibility arises from the conflicting sign structures of the representations $(+-+)$ and $(++-)$ corresponding to the functions $i(r)$ and $s(r)$, respectively.

If this property of the saving function is preserved while the corresponding property of the investment function is omitted, the computational algorithm yields only three admissible steady-state scenarios due to insufficient structural information.

Therefore, this property is incorporated implicitly through the function $q(r)=i(r)-s(r)$. The resulting scenario and transition sets, however, do not differ from those obtained for the previously considered models and therefore do not provide any new qualitative information, as mentioned above.
\end{rmk}

\begin{rmk}
The identical scenario and transition sets corresponding to the models N1.1, N1.2, and N1.3 indicate that the resulting qualitative behaviour in the considered models does not depend on the curvature of the investment and saving functions with respect to aggregate income, $i(y)$ and $s(y)$, but only on their monotonic increase. Therefore, incorporating the Kaldor conditions does not provide any additional qualitative information from the perspective of trend-based modelling.

Similarly, the implementation of the function $q(r)=i(r)-s(r)$, implicitly incorporating $s(r)$, does not generate any new qualitative behaviour. The model N1.1 and the models N2.1, N2.2, and N2.3 therefore describe the same qualitative situation. In detail, the sign-based representation of $i(r)$ is ($I_r$~D$I_r$~DD$I_r$) = $(+-+)$ and that of $s(r)$ is ($S_r$~D$S_r$~DD$S_r$) = $(++-)$. Consequently, the representation of $-s(r)$ is $(--+)$. After subtraction, the sign-based representation of $q(r)$ becomes $(+-+)$, $(0-+)$, or $(--+)$, see Table~\ref{tab:sign-based_addition}. Since trend-based modelling considers admissible scenarios and transitions without assigning them explicitly to the past, present, or future, and since the subtraction in $q(r)=i(r)-s(r)$ changes only the sign of the function value without affecting the monotonicity or curvature trends, the resulting scenario and transition sets coincide for all these models.

Therefore, the resulting scenario and transition sets exhibit qualitative robustness with respect to the considered nonlinear model representations. In other words, the qualitative behaviour obtained by trend-based modelling remains invariant under the considered nonlinear model variations.
\end{rmk}

\subsubsection{Models with an Additional Inequality}
In trend-based modelling, it is relatively simple to incorporate additional information about modelled systems by adding a qualitative inequality to the system of equations. For illustration, we add the inequality $DY<DR$ to the linear model L1 and the nonlinear model N1.1. Such an inequality may describe different adjustment speeds in the goods and money markets, where agents react faster in the money market than in the goods market.

In the model including an inequality, a scenario is regarded as admissible if the qualitative trends corresponding to the variables $y(t)$ and $r(t)$ are compatible with the inequality $DY<DR$ under admissible smooth evolution of these functions while preserving their curvature. Scenarios for which the inequality necessarily leads to an unavoidable intersection between these functions, followed by a violation of the inequality $DY<DR$, are excluded from the resulting scenario set. For instance, constant trends corresponding to a steady state may represent coinciding function values; thus, a violation of the inequality cannot be excluded from the qualitative information alone, and the corresponding scenario is rejected.

\begin{figure}[ht]
  \centering
  \includegraphics[height=6cm]{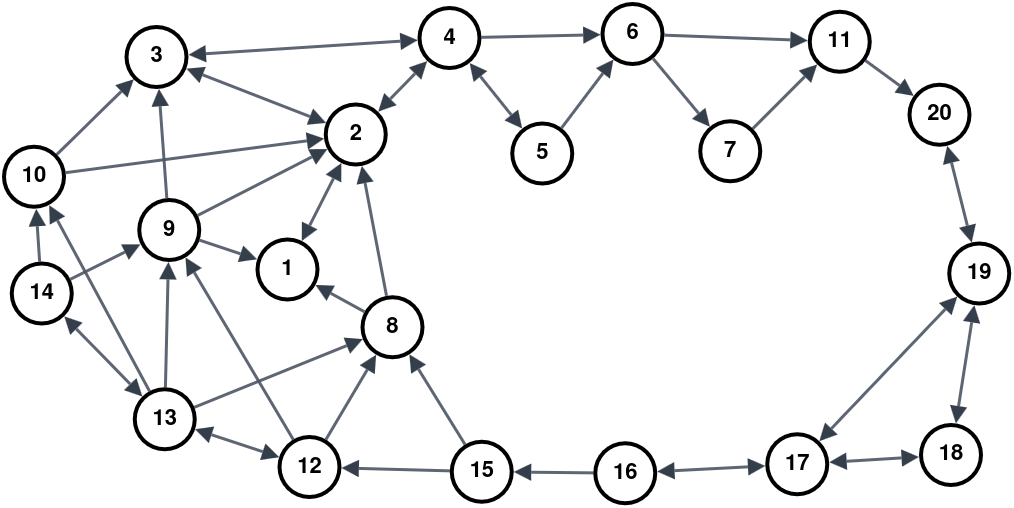}
  \caption{Transition graph corresponding to L1 with an additional inequality}
  \label{fig:transition_graph_L1_inequality}
\end{figure}

\begin{figure}[ht]
  \centering
  \includegraphics[height=5cm]{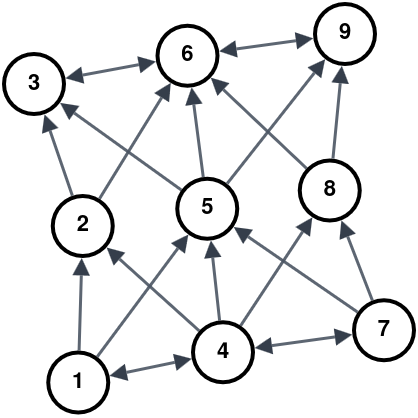}
  \caption{Transition graph corresponding to N1.1 with an additional inequality}
  \label{fig:transition_graph_N1_1_inequality}
\end{figure}

The scenario sets corresponding to L1 and N1.1 with this inequality contain 20 and 9 scenarios and are listed in Tables~\ref{tab:scenarios_L1_inequality} and~\ref{tab:scenarios_N1_1_inequality}, respectively. Naturally, the scenario sets of these models with the additional inequality form subsets of the corresponding scenario sets obtained without this inequality, see Tables~\ref{tab:scenarios_L1} and~\ref{tab:scenarios_nonlinear}. The steady-state scenarios are excluded, as explained above. The corresponding transition graphs are shown in Figures~\ref{fig:transition_graph_L1_inequality} and~\ref{fig:transition_graph_N1_1_inequality}. In Figure~\ref{fig:transition_graph_L1_inequality}, we can observe some cycles, for example, $1 \rightarrow 2 \rightarrow 3 \rightarrow 4 \rightarrow 6 \rightarrow 11 \rightarrow 20 \rightarrow 19 \rightarrow 17 \rightarrow 16 \rightarrow 15 \rightarrow 8 \rightarrow 1$, corresponding to the cycle in L1 without the additional inequality $1 \rightarrow 2 \rightarrow 3 \rightarrow 6 \rightarrow 10 \rightarrow 20 \rightarrow 33 \rightarrow 30 \rightarrow 26 \rightarrow 25 \rightarrow 24 \rightarrow 14 \rightarrow 1$, see Table~\ref{tab:scenarios_L1} and Figure~\ref{fig:transition_graph_L1}.

\begin{table}[ht]
\centering
\begin{tabular}{ccc c ccc c ccc}
\hline
No. & $Y$     & $R$     && No. & $Y$     & $R$     && No. & $Y$     & $R$ \\
\hline
1   & $+++$   & $+++$   && 8   & $+++$   & $+~0~+$ && 15  & $+~0~+$ & $+-+$ \\
2   & $++\,0$ & $+++$   && 9   & $++\,0$ & $+~0~+$ && 16  & $+-+$   & $+-+$ \\
3   & $++-$   & $+++$   && 10  & $++-$   & $+~0~+$ && 17  & $+-\,0$ & $+-+$\\
4   & $++-$   & $++\,0$ && 11  & $+--$   & $+~0~-$ && 18  & $+--$   & $+-+$ \\
5   & $++-$   & $++-$   && 12  & $+++$   & $+-+$   && 19  & $+--$   & $+-\,0$\\
6   & $+~0~-$ & $++-$   && 13  & $++\,0$ & $+-+$   && 20  & $+--$   & $+--$ \\
7   & $+--$   & $++-$   && 14  & $++-$   & $+-+$   && \\
\hline \hline
\end{tabular}
\vspace{0.2cm}
\caption{Scenario set corresponding to L1 with an additional inequality}
\label{tab:scenarios_L1_inequality}
\end{table}

\begin{table}[ht]
\centering
\begin{tabular}{ccc c ccc c ccc}
\hline
No. & $Y$ & $R$       && No. & $Y$     & $R$       && No. & $Y$     & $R$ \\
\hline
1 & $+++$ & $+-+$     && 4   & $++\,0$ & $+-+$     && 7   & $++-$ & $+-+$ \\
2 & $+++$ & $\,0\,-+$ && 5   & $++\,0$ & $\,0\,-+$ && 8   & $++-$ & $\,0\,-+$ \\
3 & $+++$ & $--+$     && 6   & $++\,0$ & $--+$     && 9   & $++-$ & $--+$ \\
\hline \hline
\end{tabular}
\vspace{0.2cm}
\caption{Scenario set corresponding to N1.1 with an additional inequality}
\label{tab:scenarios_N1_1_inequality}
\end{table}

Generally, transition graphs may contain terminal nodes or terminal components from which no directed edges lead outside the corresponding part of the graph. An example is given by the steady-state scenarios corresponding to L2, L3, and the nonlinear models. Since steady-state scenarios represent economic equilibria or stable situations in the considered models, such terminal nodes may be interpreted as economically favourable outcomes.

By contrast, the nodes no.~3, 6, and 9 corresponding to N1.1 with the additional inequality form a terminal component in the sense that transitions from other scenarios may enter this subset, while, once reached, no transitions lead outside it. These scenarios correspond to negative decreasing convex evolution of the interest rate together with positive increasing evolution of aggregate income. Depending on the economic interpretation, negative interest rates may represent either an undesirable economic situation or a possible stabilisation mechanism, while positive increasing evolution of aggregate income is generally regarded as economically favourable.

Thus, identifying such terminal structures and their associated qualitative behaviour may provide useful information about potentially favourable or unfavourable long-term developments of the modelled system.

\section{Comparative Analysis}
\label{comparative_analysis}

In Sections~\ref{phase_plane_analysis} and~\ref{trend_based_qualitative_analysis}, the well-established two-dimensional dynamic \mbox{IS--LM} macroeconomic model is qualitatively analysed in its linear and nonlinear forms using two different approaches. The first is the phase-plane approach based on a system of autonomous differential equations specified in Section \ref{system_of_differential_equation}, while the second is the trend-based approach following from the sign-based representation specified in Section~\ref{sign-based_representation}. 

In this section, we compare the information provided by these two approaches, i.e., phase-plane analysis and trend-based analysis. Both methods are qualitative, but in different senses. Phase-plane analysis works with concrete types of functions or with functions satisfying certain qualitative properties together with numerical ranges of constants and parameters, rather than only with explicit functions and fixed numerical values. From this perspective, phase-plane analysis is regarded as qualitative. Trend-based analysis operates at a different level of abstraction. This method does not work with concrete function types or numerical ranges; instead, it works only with sign triplets associated with functions, variables, and constants. Thus, trend-based analysis is regarded as qualitative in a generalized sense while using substantially more limited information.

In the previously presented analyses, we observed the key properties of these methods, which are briefly summarized in Table~\ref{tab:methods_properties}. Each property is explained in detail below, but for brief orientation we provide this table. Since both approaches are qualitative in different senses, the listed properties should be understood in a broad general sense and in the meaning described below rather than in the classical sense of dynamical systems or another well-established field.

\begin{table}[ht]
\centering
\begin{tabular}{lll}
\hline
            & Phase-plane analysis                   & Trend-based analysis \\
\hline
Method      & Information-exact                      & Information-robust \\
Abstraction & Function-level abstraction             & Sign-level abstraction \\
Inputs      & Systems of differential equations with & Sign-based representations \\
            & concrete types of functions, ranges    & of systems, variables, and \\
            & of parameters and constants            & functions \\
Outputs     & Particular types of solutions          & Scenario and transition sets \\
Results     & Local                                  & Global \\
\hline \hline
\end{tabular}
\vspace{0.2cm}
\caption{A brief overview of key phase-plane and trend-based properties}
\label{tab:methods_properties}
\end{table}

The phase-plane analysis operates with exactly specified mathematical structures, such as systems of differential equations, concrete types of functions, derivatives, Jacobians, determinants, eigenvalues, and numerical ranges of parameters and constants. We refer to this characteristic as information-exact. By contrast, the term information-robust expresses the ability of trend-based analysis to operate under substantially weaker informational assumptions. This method does not require explicit functional forms or numerical ranges and instead works only with sign triplets associated with variables, functions, and constants while still providing meaningful qualitative information about admissible scenarios and transitions.

The abstraction level reflects the range of mathematical objects with which a given method operates. Function-level abstraction refers to functions, parameters, derivatives, differential equations, singularities, and related mathematical structures. By contrast, sign-level abstraction operates only with sign triplets containing the signs of function or variable values together with the signs of their first- and second-order derivatives. This information may be incomplete while still providing meaningful qualitative results.

The inputs of the phase-plane qualitative analysis of the considered \mbox{IS--LM} models are particular functions (see the linear models) or concrete types of these functions (see the nonlinear models), together with specified ranges of parameters and constants forming the corresponding autonomous systems of differential equations. The outputs are classifications of singular points exhibited by these models. By contrast, trend-based analysis requires only sign-based representations of model elements in the form of sign triplets together with sign-based representations of the corresponding systems of differential equations as inputs. The outputs are scenario and transition sets representing admissible past, present, and future possible qualitative evolutions of the model variables.

The results of phase-plane analysis are strictly local, especially in nonlinear models, due to the local character of the dynamics in the neighbourhood of singular points. In particular, the linearised system of differential equations only approximates the behaviour in a neighbourhood of a singular point. By contrast, the global character of trend-based analysis is understood in the sense that it generates sets of admissible scenarios together with transitions between them, representing possible qualitative evolutions of the considered system. This method does not identify the number of singular points, nor provide their dynamical classification. However, it enables the observation of the existence of steady states corresponding to singular points together with admissible qualitative time evolutions described by sign triplets expressing the value, monotonicity, and curvature of the observed variables. Precisely, phase-plane analysis provides a local dynamical classification of singular points, distinguishing particular types of local behaviour such as sinks, saddles, or sources. By contrast, trend-based analysis aggregates these local dynamical possibilities into admissible steady-state configurations enriched by transition structures. Consequently, a single steady-state scenario in trend-based analysis may correspond to several qualitatively different local dynamical structures identified by phase-plane analysis. In this sense, trend-based analysis captures a more global qualitative structure of the admissible dynamics rather than particular local classifications around singular points.

\subsection{Comparison in the Linear Case}
Phase-plane analysis of the linear \mbox{IS--LM} model works with the system of differential equations \eqref{linear_IS-LM_model}, which already includes the particular specification of the model functions given in \eqref{linear_relations}. Moreover, the corresponding functional forms are naturally incorporated into the model through additive relations preserving the linear character of the system. The classification of singular points associated with the linear model is provided in Tables~\ref{tab:hyperbolic_eq_points} and~\ref{tab:non-hyperbolic_eq_points}. In fact, except for the degenerate case of overlapping IS and LM curves, one linear system of differential equations may lead to one singular point of different dynamical types depending on the numerical ranges of the model parameters.

Besides this, for computational and interpretative reasons, trend-based analysis in the linear case works with four systems of equations, L1, L2, L3, and L4, expressed in the sign-based representation, see Table~\ref{tab:systems_representations}. The systems L2, L3, and L4 together form the system L1. In fact, the models L2 and L3 correspond to the situations illustrated in the first and second rows of Table~\ref{tab:hyperbolic_eq_points}, respectively. In the transition graphs corresponding to L2 and L3, see Figures~\ref{fig:transition_graph_L2} and~\ref{fig:transition_graph_L3}, we can observe trend-based qualitative structures resembling a sink, i.e., a stable node or a stable focus. Precisely, the steady-state scenario can be reached from each node without passing through any cycle, which qualitatively resembles trajectories associated with a stable node. Similarly, the steady-state scenario can be reached from each node after passing through one or more cycles beforehand, which qualitatively resembles oscillatory trajectories associated with a stable focus. Since trend-based analysis works only with positivity, monotonicity, and curvature, oscillatory behaviour with both constant and varying amplitudes may be represented by cycles in the transition graph. Thus, in summary, both types of trend-based qualitative behaviour are aggregated into a single transition graph, although they correspond to two different types of singular points represented by different phase portraits in phase-plane analysis. 

The system L4 corresponds to the cases described by the third or fourth row of Table~\ref{tab:hyperbolic_eq_points}, or the first row of Table~\ref{tab:non-hyperbolic_eq_points}, or the first row of Table~\ref{tab:overlapping_curves_eq_points}. Since this model incorporates several possible types of singular points, the corresponding transition graph is more complicated, see Figure~\ref{fig:transition_graph_L4}. In this transition graph, we can observe trend-based qualitative structures resembling qualitative behaviours associated with all singular point types identified by phase-plane analysis, i.e., sink, source, saddle, centre, and attracting or repelling sets in the degenerate case. In particular, there are cycles not containing a steady-state scenario, qualitatively resembling periodic behaviour associated with a centre, and sequences of nodes and directed edges qualitatively resembling sink-like and source-like behaviour in the sense described above. Consequently, we can observe here the property that the steady-state scenario can be reached from two scenarios and can be left through two different scenarios, with corresponding sign triplets qualitatively resembling trajectories associated with a saddle. Finally, the degenerate case of overlapping IS and LM curves is characterised by infinitely many attracting or repelling singular points lying on a straight line, namely on the IS or LM curve. Since these singular points attain positive values, this straight line is represented in trend-based analysis by a single steady-state scenario $\{(+~0~0),(+~0~0)\}$. The attracting or repelling character is represented by directed edges leading to or from this node. Since trend-based analysis admits only one steady-state scenario $\{(+~0~0),(+~0~0)\}$ in these linear models, all described qualitative behaviours primarily identified by phase-plane analysis are aggregated around one node associated with this steady-state configuration from the perspective of trend-based analysis. Therefore, this transition graph allows transitions between various types of trend-based qualitative behaviour, for example from sink-like to saddle-like behaviour. This phenomenon may be interpreted as capturing changes in the numerical ranges of model parameters leading to changes in the corresponding types of singular points.

The system L1 corresponds to all cases described in Table~\ref{tab:hyperbolic_eq_points}, or the first row of Table~\ref{tab:non-hyperbolic_eq_points}, or the first row of Table~\ref{tab:overlapping_curves_eq_points}. As can be seen in Figure~\ref{fig:transition_graph_L1}, the corresponding transition graph is not sufficiently transparent; however, the previously described trend-based qualitative structures resembling all mentioned singular point types can be found here.

Thus, we can see that in these linear cases trend-based analysis provides a more general qualitative representation integrating the considered possibilities into one transition graph, while phase-plane analysis provides a more detailed description decomposed into several phase portraits.

As mentioned in Remark \ref{rmk:linear_model}, the linear model can otherwise be expressed using the sign-based representation of the general system N1 together with the function property representations $(++0)$ and $(+-0)$ for increasing and decreasing linear functions, respectively. Such a sign-based representation does not lead to sufficiently interpretative results. Only the steady-state scenario is admissible in this case. The general formulation admits a wide range of possible behaviours, but due to insufficient structural information, most dynamic scenarios cannot be either excluded or confirmed. Therefore, the trend-based modelling algorithm is not able to determine them uniquely and retains only the admissible scenario corresponding to the steady state. This steady-state scenario remains compatible with the available qualitative information and does not require any additional structural assumptions. Thus, the general sign-based formulation is not poorer in the sense of possible behaviours, but rather less decidable from the viewpoint of trend-based analysis. The resulting steady-state scenario indicates that the existence of an equilibrium is consistent with the available information, although no further qualitative properties can be unambiguously determined. On the other hand, the phase-plane analysis provides a sufficiently clear description of model dynamics in this case. Since this analysis uses only first-order partial derivatives of the model functions, the classification of singular point types through the Jacobian matrix can still be reliably performed in the linear model expressed by this general specification \eqref{IS-LM_model} using the properties \eqref{economic_is_y}, \eqref{economic_is_r}, and \eqref{economic_l_yr}. In fact, the classification of a singular point agrees with that listed in Tables~\ref{tab:hyperbolic_eq_points} and~\ref{tab:non-hyperbolic_eq_points}, while the degenerate case of overlapping IS and LM curves corresponds to the situation described in Table~\ref{tab:overlapping_curves_eq_points}.

\subsection{Comparison in the Nonlinear Case}
Perhaps even more illustrative models for demonstrating differences between these qualitative approaches are the nonlinear \mbox{IS--LM} models, especially the Kaldor-type \mbox{IS--LM} model. From the perspective of phase-plane analysis, three singular points may emerge, where the first and third points are sinks and the middle point is a saddle, see Table~\ref{tab:Kaldor_IS-LM_model_eq_points}. By contrast, from the perspective of trend-based analysis, resulting in the scenario set listed in Table~\ref{tab:scenarios_nonlinear} and the transition graph illustrated in Figure~\ref{fig:transition_graph_nonlinear}, only steady-state scenarios associated with attractor-like qualitative behaviour are observed. Precisely, the scenarios no.~10, 11, and 12 are steady-state scenarios differing only in the sign appearing in the first position of the triplets corresponding to the interest rate $r$ and exhibiting such a behaviour with directed edges leading only towards these steady states, either through non-cyclic paths or after passing through degenerate cycles beforehand. As a degenerate cycle, we consider cycles such as $1 \rightarrow 4 \rightarrow 7 \rightarrow 4 \rightarrow 1$ or $3 \rightarrow 6 \rightarrow 9 \rightarrow 6 \rightarrow 3$ in this transition graph. Such a degenerate cycle followed by a trend path directed towards a steady-state scenario may be interpreted as a trend-based object exhibiting monotone trend evolution in the variable $r$ and oscillatory trend evolution with curvature changes in the variable $y$. This behaviour does not appear to have a direct correspondence with the classical classification of singular points.

The absence of the intermediate saddle-like structure in trend-based analysis is primarily caused by the inability of the current computational algorithm to capture changes in curvature within a single function. In the Kaldor-type \mbox{IS--LM} model, the intermediate equilibrium structure arises from the S-shaped form of the IS curve, where convex and concave regions coexist. This S-shaped IS curve follows from the curvature change in the model functions $i(y)$ and $s(y)$, see \eqref{Kaldor_is}. Unfortunately, the computational algorithm allows only one qualitative sign-triplet structure for the corresponding model functions and therefore does not allow the two different triplets required by \eqref{Kaldor_i} and \eqref{Kaldor_s}. Consequently, this internal geometric change is not preserved within the trend-based abstraction.

Therefore, the nonlinear trend-based representation effectively reduces to the qualitative structure corresponding to decreasing IS and increasing LM curves, which primarily leads to attractor-like behaviour. This corresponds to the situations described in the first and third columns of Table~\ref{tab:Kaldor_IS-LM_model_eq_points}. The three configurations of steady-state scenarios differ only in the sign of the interest rate value. They may be visualised by relative vertical shifts of the IS and LM curves leading to equilibrium points located in regions with positive, zero, or negative values of $r$. However, in the Kaldor-type \mbox{IS--LM} model, the regions with decreasing IS curves and increasing LM curves prevail; thus, this type of dynamic behaviour may be regarded as dominant in such a model.

\newpage
Generally, it seems that trend-based analysis in all considered nonlinear \mbox{IS--LM} models reduces the sets of admissible scenarios and admissible transitions between them to the case corresponding to attractor-like behaviour, since these resulting sets are identical for all these models. This situation agrees with the standard monotonically decreasing IS curve and monotonically increasing LM curve corresponding to the standard qualitative shapes of these curves from the economic point of view. Non-standard situations exhibiting different qualitative shapes of IS and LM curves, such as those describing a liquidity trap, are not captured by the trend-based abstraction following from the considered sign-based representations of the nonlinear models. By contrast, phase-plane analysis provides in these cases a relatively detailed description of local dynamics for all possible qualitative shapes of IS and LM curves. Nevertheless, some extension of the computational algorithm may ensure the possibility of working with more complicated sign-based representations of model functions, and thus richer trend-based dynamics may be obtained.

\section{Conclusion}

Two qualitative analytic methods based on different forms of abstraction are applied to several versions of a well-established dynamic macroeconomic model. Both methods, phase-plane analysis and trend-based analysis, are qualitative in character but in different senses. The well-established phase-plane analysis provides sufficiently clear results for local dynamics with the use of hard data. By contrast, trend-based analysis also works in situations with incomplete information and aggregates local results into one illustrative representation in the form of a transition graph, but is limited by the current computational algorithm.

Phase-plane analysis provides relatively precise results in the form of classifications of singular points. However, the economic interpretation of such results may be difficult for non-mathematicians. On the other hand, the illustration of trend-based analysis results in the form of transition graphs is more accessible for non-experts. The meaning of terminal nodes or favourable outcomes is sufficiently clear even for an economic audience. Similarly, the formulation of inputs required in phase-plane analysis in the form of functional relations may be less natural for non-mathematicians. By contrast, the formulation of inputs required in trend-based analysis in the form of vague heuristics, such as a decelerating increase of some economic quantity, may be more desirable.

After a detailed analysis and comparison of these two qualitative methods applied to two-dimensional dynamic \mbox{IS--LM} models, a suitable combination of these qualitative approaches may provide a meaningful and broader picture of model dynamics. As we can observe, in the linear models, the resulting trend-based structures exhibit a close correspondence with the classical classification of singular points obtained by phase-plane analysis. In the nonlinear models, however, this correspondence becomes less direct. While attractor-like qualitative behaviour can still be identified, the resulting transition structures may contain admissible trend configurations that are not naturally described by the local classification of singular points. This observation suggests that trend-based analysis may capture additional global qualitative structures whose relationship to classical phase-plane concepts remains an open question.

\section*{Acknowledgements}
The research was supported by the Ministry of Education, Youth and Sports of the Czech Republic (MSMT CR) under RVO funding for IC47813059 and IC00216305.

\section*{Statement}
During the preparation of this work the authors used ChatGPT (OpenAI) to improve the clarity of the English language. After using this tool, the authors reviewed and edited the content as needed and take full responsibility for the content of the publication.

\bibliographystyle{plain}
\bibliography{references}  

\end{document}